\documentclass[10pt]{article}         
\usepackage{amsmath,amsfonts,amssymb,mathrsfs,amsthm,url}
\usepackage[utf8]{inputenc}
\usepackage[T1]{fontenc}
\usepackage{lmodern}
\usepackage{hyperref}
\usepackage{cite}

\newtheorem{thm}{Theorem}[section]
\newtheorem{THM}{Theorem}

\theoremstyle{definition}
\newtheorem{defn}[thm]{Definition}
\newtheorem{prop}[thm]{Proposition}
\newtheorem{lem}[thm]{Lemma}
\newtheorem{cor}[thm]{Corollary}

\theoremstyle{remark}
\newtheorem{rem}[thm]{Remark}

\newtheorem{eg}[thm]{Example}

\title{Open Homomorphisms between $m$-step Solvable Galois Groups of Number Fields}
\author{Yu Mao \\
\texttt{China Three Gorges University}\\
\texttt{maoyu@ctgu.edu.cn}}

\date{}
\begin{document}
\maketitle
\begin{abstract}
Let $K,L$ be number fields. In this paper, we prove that, for $m\geq 2$, a homomorphism $\phi_{m+2}:G_K^{m+2}\to G_L^{m+2}$ whose depth-two projection has open image gives rise to a unique equivariant field embedding $L_m\hookrightarrow K_m$. We first establish a sharper version of the $m$-step solvable cyclotomic compatibility criterion of Mao and Sa\"idi. This criterion also holds for $m=1$, with uniqueness on the base field. We then use Frobenius detection for additive $\ell$-adic characters and an auxiliary character module to show that a homomorphism $G_K^4\to G_L^4$ with open image at depth two is automatically compatible with the cyclotomic characters.
\end{abstract}
\newpage
\tableofcontents
\section*{Conventions}
Throughout this paper, we shall use the following conventions.
\begin{itemize}
    \item For a field $k$, write $\bar{k}$ for its separable closure contained in some algebraic closure $k^{\text{alg}}$ of $k$.
    \item For a profinite group $G$, put $G^{[0]}=G$ and $G^{[i]}=\overline{[G^{[i-1]},G^{[i-1]}]}$ for $i\geq1$. We call $G^i=G/G^{[i]}$ the maximal $i$-step solvable quotient of $G$.
    \item For an abelian group $G$, and an integer $ n \geq 2$, we shall write
    $$
    G[n]:= \text{ker}(G \xrightarrow{\times n} G)
    $$
    if $G$ is viewed additively.
    \item For a field $k$, write $G_k^m$ for the maximal $m$-step solvable quotient of $G_k$ for some $m \geq 0$. We write $k_m=\bar{k}^{G_k^{[m]}}$ for the maximal $m$-step solvable extension of $k$. In particular, $k_m = k_{m-1}^{\text{ab}}$ for all $m \geq 1$ (where $k_0 := k$).
    \item We write $\mathfrak{Primes}$ for the set of all rational primes. 
    \item For an algebraic extension $k/\mathbb Q$, write $\mathscr{P}_k$ for the set of all primes of $k$. Moreover, write $\mathscr{P}_k^{<\infty}$ for the set of all non-archimedean primes of $k$ and $\mathscr{P}_k^{\infty}$ for the set of all archimedean primes of $k$. 
    \item Let $k$ be a number field and let $v\in\mathscr P_k$ be a prime. We write $k_v$ for the completion of $k$ at $v$. 
    \item Let $k$ be a number field, and let $L/k$ be a Galois extension, let $v \in \mathscr{P}_k^{<\infty}$ and $w \in \mathscr{P}_L^{<\infty}$ lying above $v$. We write $D_w(L/k)$ (resp. $I_w(L/k)$) for the decomposition subgroup of $\text{Gal}(L/k)$ at $w$ (resp. the inertia subgroup of $\text{Gal}(L/k)$ at $w$). If there is no ambiguity, we shall simply write $D_w$ (resp. $I_w$) for this decomposition group (or even $D_v$ if we only need the conjugacy class).
    \item For a number field or a $p$-adic local field $k$, write $\chi_k$ for the cyclotomic character $\chi_k: G_k \to \widehat{\mathbb{Z}}^\times$. For some prime number $\ell$, write $\chi_{k,\ell}$ for the $\ell$-adic cyclotomic character $\chi_{k,\ell}:G_k \to \mathbb{Z}_\ell^\times$.
    \item All homomorphisms between profinite groups in this paper are assumed to be continuous.
    \item We write $P=G_{\mathbb{Q}}$ and write $\mathbb{Q}^{(j)}$ for the field denoted by $k_j$ above when $k=\mathbb{Q}$. The superscript distinguishes the solvable layer $\mathbb{Q}^{(p)}$ from the local field $\mathbb{Q}_p$ when the index is a prime $p$. Thus $P^j=\operatorname{Gal}(\mathbb{Q}^{(j)}/\mathbb{Q})$. If $\phi_n:G_K^n\to G_L^n$ and $1\leq j\leq n$, we write $\phi_j$ for the induced map at depth $j$.
\end{itemize} 
\section{Introduction}
In \cite{Uchida3}, Uchida conjectured that every open homomorphism between solvably closed Galois groups of number fields arises uniquely from a field embedding. This is the hom-form of the birational anabelian conjecture for number fields, whose isom-form is the Neukirch--Uchida theorem. The following is the main theorem of \cite{OAI}.
\begin{thm}[cf. Theorem 1.1 in \cite{OAI}]
Let $F_1,F_2$ be number fields, and let $E_i/F_i$ be solvably closed Galois extensions for $i \in \{1,2\}$. Then every open homomorphism
$$
\alpha: \text{Gal}(E_1/F_1) \to \text{Gal}(E_2/F_2)
$$
arises uniquely from a field embedding $\beta: E_2 \hookrightarrow E_1$ such that $g \circ \beta = \beta \circ \alpha(g)$ for all $g \in \text{Gal}(E_1/F_1)$.
\end{thm}
The proof of \cite{OAI} uses Frobenius detection for additive $\ell$-adic characters and deduces full cyclotomic compatibility from logarithmic compatibility. On the isomorphism side, Sa\"idi and Tamagawa proved the following finite-step solvable version of the Neukirch--Uchida theorem.
\begin{thm}[cf. Theorem 2 in \cite{ST1}]
Let $m \geq 0$ be an integer and let $K,L$ be number fields. Let $\tau_{m+3}: G_K^{m+3} \xrightarrow{\sim} G_L^{m+3}$ be an isomorphism. \par 
(i) There exists an isomorphism $\sigma_m:K_m \xrightarrow{\sim} L_m$ such that $\tau_m(g) = \sigma_mg\sigma_m^{-1}$ for all $g \in G_K^m$ where $\tau_m:G_K^m \xrightarrow{\sim} G_L^m$ is the isomorphism induced by $\tau_{m+3}$. In particular, $\sigma_m$ restricts to an isomorphism $\sigma:K \xrightarrow{\sim} L$. \par 
(ii) Assume that $m \geq 2$ (resp. $m = 1$). Then the isomorphism $\sigma_m:K_m \xrightarrow{\sim} L_m$ (resp. $\sigma:K \xrightarrow{\sim} L$ induced by $\sigma_1:K_1 \xrightarrow{\sim} L_1$) in (i) is uniquely determined by the property that $\tau_m(g) = \sigma_mg\sigma_m^{-1}$ for every $g \in G_K^m$.
\end{thm} 
It is natural to ask for a finite-step solvable version of Theorem 1.1. Corato and Sa\"idi obtained such a result under local and image conditions at input depth $m+3$ in \cite{CS}. Our main result removes those conditions in favour of openness at depth two and reduces the input depth to $m+2$.
\begin{THM}\label{intro:main}
Let $m\geq 2$ and let $K,L$ be number fields. Let $\phi_{m+2}:G_K^{m+2}\to G_L^{m+2}$ be a homomorphism. Then the following conditions are equivalent:
\begin{enumerate}
\item $\phi_2(G_K^2)$ is open in $G_L^2$;
\item $\chi_L\circ\phi_{m+2}=\chi_K$;
\item there exists a field embedding $\tau_m:L_m\hookrightarrow K_m$ satisfying $g\circ\tau_m=\tau_m\circ\phi_m(g)$ for every $g\in G_K^m$.
\end{enumerate}
The embedding in (3) is unique and sends $L$ into $K$.
\end{THM}
\begin{rem}
The openness condition in Theorem~\ref{intro:main} is strictly weaker than openness of the given lift $\phi_{m+2}$. Indeed, Example~\ref{eg:nonopen-lift} gives, for every $m\geq2$, an injective endomorphism of $G_{\mathbb Q}^{m+2}$ which induces the identity on $G_{\mathbb Q}^{m+1}$ and has infinite-index image.
\end{rem}
The first part of our proof gives a sharper version of the cyclotomic criterion in \cite[Theorem 1.5]{MS}: the input depth is $m+2$ rather than $m+3$, and the conditional criterion does not require openness.
\begin{THM}\label{intro:criterion}
Let $m\geq 1$ and let $\phi_{m+2}:G_K^{m+2}\to G_L^{m+2}$ be a homomorphism. Then $\chi_L\circ\phi_{m+2}=\chi_K$ if and only if there exists an embedding $\tau_m:L_m\hookrightarrow K_m$ satisfying $g\circ\tau_m=\tau_m\circ\phi_m(g)$ for all $g\in G_K^m$. In this case $\phi_m$ is open and the restriction $\tau_m|_L:L\hookrightarrow K$ is unique. If $m\geq2$, the embedding $\tau_m$ is unique. If $m=1$, the equivariant embeddings form a torsor under $G_L^{\rm ab}$ by precomposition.
\end{THM}
An application of Theorem~\ref{intro:criterion} to isomorphisms gives a two-step loss isom-form, including an isomorphism at the abelian layer from input depth three (Corollary~\ref{cor:isom}).
\begin{THM}\label{intro:automatic}
Every homomorphism $\beta:G_K\to P^4$ whose projection to $P^2$ has open image satisfies $\chi_{\mathbb Q}\circ\beta=\chi_K$. Consequently a homomorphism $G_K^4\to G_L^4$ whose projection at depth two has open image is compatible with the cyclotomic characters.
\end{THM}
Theorems~\ref{intro:criterion} and \ref{intro:automatic} together imply Theorem~\ref{intro:main}. For isomorphisms at input depth three, full cyclotomic compatibility is already available by \cite[Theorem 1.26]{ST1}, so Theorem~\ref{intro:criterion} applies at $m=1$. For homomorphisms the additional arithmetic input is Theorem~\ref{intro:automatic}, which uses depth four.

The argument for Theorem~\ref{intro:automatic} follows the passage from logarithmic to full cyclotomic compatibility in \cite{OAI}. Its extension to the normaliser of the kernel \cite[Lemma 4.1]{OAI}, however, uses openness of the image in the full target and the realisation of automorphisms of that open subgroup by ambient conjugation. Neither input is available in this generality at finite depth: Example~\ref{eg:nonopen-lift} shows that the image of the given lift need not be open, and Remark~\ref{rem:normaliser-obstruction} shows that even the whole finite-step target can have automorphisms which are not realised by ambient conjugation. Instead, in Section~4 we use the auxiliary character module of \cite[Lemma 1.7]{Tsuji}. Local decomposition-group containment and Frobenius detection give injectivity of pullback on the primitive summand needed in the proof, without requiring an open image at depth three or an extension to a normaliser.

This paper is organised as follows. In Section~2, we analyse the images of pro-$\ell$ Sylow subgroups of decomposition groups. In Section~3, we prove the cyclotomic criterion. In Section~4, we establish logarithmic cyclotomic compatibility at input depth four. In Section~5, we deduce full cyclotomic compatibility and complete the proof of the main theorem. Appendix~\ref{app:depth-three} records a logarithmic consequence at input depth three.

\section{Local Analysis of Decomposition Groups}
We shall briefly review the parts of the local theory from \cite{ST1} and \cite{CS} we need. \par 
Let $K$ be a number field, and let $\ell \in \mathfrak{Primes}$ be a prime number. Let $j \geq 2$ be an integer and let $\mathfrak{p} \in \mathscr{P}_K^{<\infty}$ be a prime not lying above $\ell$ and choose prolongations $\mathfrak{p}_j \in \mathscr{P}_{K_j}^{<\infty},\bar{\mathfrak{p}} \in \mathscr{P}_{\bar{K}}^{<\infty}$ of $\mathfrak{p}$. Then Proposition 1.1 in \cite{ST1} implies that the natural surjection
$$
D_{\bar{\mathfrak{p}}}(\bar{K}/K) \twoheadrightarrow D_{\mathfrak{p}_j}(K_j/K)
$$
preserves the maximal tame quotient, so does the maximal unramified quotient. Furthermore, a pro-$\ell$ Sylow $S \subset D_{\mathfrak{p}_j}$ is mapped isomorphically onto a pro-$\ell$ Sylow subgroup of the maximal tame quotient of $D_{\mathfrak{p}_j}$. In particular, we have an exact sequence
$$
1 \to J \to S \to S^{\text{unr}} \to 1
$$
where both $J$ and $S^{\text{unr}}$ are free pro-$\ell$ groups of rank $1$. Furthermore, the centraliser of $J$ in $S$ is $J$ itself. In particular, if a homomorphism is injective on $J$, it is injective on $S$.
\begin{lem}\label{lem:local}
Let $n\geq3$ and let $\alpha:G_K\to G_L^n$ be a homomorphism, with projections $\alpha_j:G_K\to G_L^j$. Let $\mathfrak p$ be a finite prime of $K$ away from $\ell$, choose a prolongation $\mathfrak P$ to $\bar K$, and let $S\subset D_{\mathfrak P}$ be a pro-$\ell$ Sylow subgroup. Assume that $\alpha_{n-1}|_S$ is injective.
\begin{enumerate}
\item There is a prime $\mathfrak q_{n-1}$ of $L_{n-1}$ away from $\ell$ such that, writing $\mathfrak q_{n-2}=\mathfrak q_{n-1}|_{L_{n-2}}$, one has
\[
 \alpha_{n-1}(S)\subset D_{\mathfrak q_{n-1}}(L_{n-1}/L),
 \qquad
 \alpha_{n-2}(D_{\mathfrak P})\subset D_{\mathfrak q_{n-2}}(L_{n-2}/L).
\]
The lower prime is independent of the choice of $S$ for the given $\mathfrak P$.
\item If $n\geq4$, the image of $J=S\cap I_{\mathfrak P}$ under $\alpha_{n-2}$ is infinite and is contained in $I_{\mathfrak q_{n-2}}(L_{n-2}/L)$.
\end{enumerate}
\end{lem}
\begin{proof}
Put $j=n-1$. The image of $S$ in $G_K^j$ is a $(\star_\ell)$-subgroup by \cite[Proposition 1.22]{ST1}. Since $\alpha_j$ is injective on this subgroup and lifts to depth $j+1$, its image satisfies $(\star_\ell)$ by \cite[Proposition 2.5]{CS}. Indeed, the isomorphism between this subgroup and its image identifies their degree-two cohomology with coefficients in $\mathbb F_\ell$. The map between their inverse images at depth $j+1$ makes the inflation square commute. Thus a class with nonzero inflation on the source has nonzero inflation on the target; no surjectivity on the inverse images is needed. Local recognition by \cite[Proposition 1.22]{ST1} then gives the upper containment, and moreover its image is open in a pro-$\ell$ Sylow subgroup of the corresponding decomposition group.

For $g\in D_{\mathfrak P}$, the images of $S$ and $gSg^{-1}$ in $G_K^j$ determine the same class in \cite[Proposition 1.23]{ST1}. Their images in $G_L^j$ also determine the same class by \cite[Proposition 2.9]{CS}. The equivariance of the class-to-prime map gives that $\alpha_{j-1}(g)$ fixes $\mathfrak q_{j-1}$. The same argument for two pro-$\ell$ Sylow subgroups proves independence. If $j-1\geq2$, projection of the target decomposition group preserves its pro-$\ell$ Sylow subgroups by \cite[Proposition 1.1(vii)]{ST1}. In an open subgroup of such a Sylow subgroup the intersection with tame inertia is the maximal normal abelian subgroup. Comparing it with $J$ gives the last assertion.
\end{proof}
We recall the definition of active primes from \cite{OAI}.
\begin{defn}
For an odd prime $\ell$ and a homomorphism $\beta:G_K\to P^4$, a prime $\mathfrak p\nmid\ell$ of $K$ is \emph{active at $\ell$} if $\beta_2(J)\neq1$, where $J=S\cap I_{\mathfrak P}$ for a pro-$\ell$ Sylow subgroup $S$ of $D_{\mathfrak P}$. This is independent of the choices by conjugacy. By \cite[Proposition 1.9]{CS}, $\beta_2(J)$ is infinite, so $\beta_2$ is injective on $J$ and on $S$ by the centraliser argument above. Lemma~\ref{lem:local} therefore applies. We call the rational prime below the lower target prime the \emph{target} of $\mathfrak p$; it is independent of the prolongation.
\end{defn}
Next we shall obtain the residue exponent comparison in a similar way as in Lemma 2.6 in \cite{OAI}.
\begin{lem}
Let $v\in \mathscr{P}_K^{<\infty}$ be active at $\ell$ for $\beta:G_K\to P^4$, with target prime $q$, and choose compatible prolongations. For $g\in D_v$, let $b(g) \in \widehat{\mathbb{Z}}$ be the residue exponent of $g$ at $v$ and let $a(g) \in \widehat{\mathbb{Z}}$ be the residue exponent of $\beta_2(g)$ at the chosen target prime of $\mathbb{Q}^{(2)}$. Then
$$
(Nv)^{b(g)} = q^{a(g)}~\text{in}~\mathbb{Z}_\ell^\times.
$$
In particular, $\chi_{\mathbb{Q},\ell}(\beta(g)) = \chi_{K,\ell}(g)$ on $D_v.$
\end{lem}
\begin{proof}
Choose a topological generator $\tau \in J$. Since $v\nmid\ell$, we have $J \xrightarrow{\sim} \mathbb{Z}_\ell(1)$. It follows from Theorem 7.5.3 in \cite{NSW} that $g\tau g^{-1}(\tau^{(Nv)^{b(g)}})^{-1}$ belongs to the wild inertia of $D_v$, a pro-$p$ group where $p:=\operatorname{char}(v)$. In the target tame quotient, the image of $\tau$ is non-trivial in the pro-$\ell$ tame inertia factor by Lemma 2.1. Write $t$ for this image. The image of the wild error term lies in the pro-$\ell$ tame inertia factor, because both $\tau$ and its conjugate map there. It is also an image of a pro-$p$ element, so it is trivial, as $p\neq\ell$. The target tame relation now gives 
$$
t^{q^{a(g)} - (Nv)^{b(g)}} = 1
$$
which gives the desired equality. The cyclotomic assertion follows since both source and target primes have residue characteristic different from $\ell$.
\end{proof}
\section{The \texorpdfstring{$m$}{m}-step Solvable Version of Hoshi's Cyclotomic Criterion}
In this section, we give the proof of Theorem 2. We first recall some definitions given in \cite{Ho3}.
\begin{itemize}
\item For a number field $T$, put $A_\ell(T)=G_T^{\rm ab}/\ell G_T^{\rm ab}$.
\item For a finite prime $v$ of $T$, write $\operatorname{FL}_\ell(v)\subset A_\ell(T)$ for the set of images of its Frobenius liftings. This set is independent of the prolongation.
\item A homomorphism $A_\ell(F)\to A_\ell(E)$ is Frobenius-preserving if, for every absolute residue degree-one prime $v$ above a cofinite set of rational primes, there is a prime $w$ of $E$ above the same rational prime and of absolute residue degree one, such that the image of $\operatorname{FL}_\ell(v)$ is contained in $\operatorname{FL}_\ell(w)$.
\end{itemize}
Fix $m\geq1$ and $\phi_{m+2}:G_K^{m+2}\to G_L^{m+2}$ such that $\chi_L\circ\phi_{m+2}=\chi_K$. Let $E/L$ run through the cofinal family of finite Galois subextensions of $L_m/L$ containing $i=\sqrt{-1}$. Let $G_F$ be the inverse image in $G_K$ of $\operatorname{Gal}(L_m/E)$ under the map induced by $\phi_m$. Then $F/K$ is finite Galois and $F\subset K_m$. Since $G_L^{[m]}\subset G_E$, we have $G_L^{[m+1]}\subset G_E^{[1]}$, hence $E^{\rm ab}\subset L_{m+1}$. Thus restriction gives a map $G_F\to\operatorname{Gal}(L_{m+1}/E)\to G_E^{\rm ab}$, and therefore a homomorphism
\[
 \theta_E:A_2(F)\longrightarrow A_2(E).
\]
\begin{lem}\label{lem:frob}
The map $\theta_E$ is Frobenius-preserving.
\end{lem}
\begin{proof}
Let $u$ be a degree-one prime of $F$ above an odd rational prime $p$ unramified in $F/\mathbb Q$. Its restriction $v$ to $K$ satisfies $F_u=K_v=\mathbb Q_p$, so we may identify their decomposition groups. Let $S$ be a pro-$2$ Sylow subgroup and $J=S\cap I_u$. By \cite[Lemma 3.1(i)]{Ho3}, $\chi_F(J)\neq1$ (already in its $p$-component), so $\phi_{m+1}(J)\neq1$. This image lies in $\operatorname{Gal}(L_{m+1}/E)$, which has no involutions: by \cite[Proposition 1.9]{CS}, an involution at depth $m+1$ projects to complex conjugation at depth $m$ and moves $i\in E$. It follows that $\phi_{m+1}|_J$ is injective, and hence $\phi_{m+1}|_S$ is injective by the centraliser argument preceding Lemma~\ref{lem:local}.

Apply Lemma~\ref{lem:local}(1) at depth $m+2$. At the upper prime $\mathfrak q_{m+1}$ we have the containment of $\phi_{m+1}(S)$, and at $\mathfrak q_m=\mathfrak q_{m+1}|_{L_m}$ the containment of $\phi_m(D_v)$. Put $w=\mathfrak q_m|_E$. Since the cyclotomic character factors through depth one, the lower containment and cyclotomic compatibility give $\chi_F(D_u)\subset\chi_E(D_w)$. By \cite[Lemma 3.1(ii),(iii)]{Ho3}, the residue characteristic of $w$ is $p$ and its absolute residue degree is one: the first image contains $\mathbb Z_p^2$, and its prime-to-$p$ projection contains $p$.

Let $g\in D_u$ be a Frobenius lifting. The image of $D_u$ in $A_2(F)\times\mathbb Z/2\mathbb Z$, with residue exponent modulo $2$ in the second coordinate, is a pro-$2$ group; it is therefore the image of $S$. Choose $s\in S$ with the same image as $g$. The upper containment gives a lift $t\in D_w\subset G_E$ of $\phi_{m+1}(s)$. Then $[t]=\theta_E([g])$ and the prime-to-$p$ part of cyclotomic compatibility gives
\[
 p^{a_w(t)}=p^{a_u(s)}\quad\text{in }\prod_{q\neq p}\mathbb Z_q^\times,
\]
where $a_w,a_u$ denote residue exponents. The map $a\mapsto p^a$ from $\widehat{\mathbb Z}$ to this product is injective: reduction modulo $p^n-1$ detects $a\bmod n$, since $p$ has order $n$ modulo $p^n-1$ for $n\geq2$. Hence $a=a_w(t)$ is odd. Choose $z\in D_w$ with residue exponent $(1-a)/2$. Then $tz^2$ is a Frobenius lifting with class $[t]$ in $A_2(E)$. This proves the required inclusion of Frobenius lifting sets.
\end{proof}

\begin{thm}[Theorem~\ref{intro:criterion}]\label{thm:criterion}
Let $m\geq1$ and let $\phi_{m+2}:G_K^{m+2}\to G_L^{m+2}$ be a homomorphism. The equality $\chi_L\circ\phi_{m+2}=\chi_K$ is equivalent to the existence of a $\phi_m$-equivariant embedding $\tau_m:L_m\hookrightarrow K_m$. Whenever the embedding exists, $\phi_m$ is open and the restriction to $L$ is unique. The embedding itself is unique for $m\geq2$; for $m=1$ the embeddings form a $G_L^{\rm ab}$-torsor.
\end{thm}
\begin{proof}
Assume first that $\chi_L\circ\phi_{m+2}=\chi_K$. Lemma~\ref{lem:frob} and \cite[Corollary 2.8]{Ho3} give a unique embedding $\tau_E:E\hookrightarrow F$ inducing $\theta_E$. For $g\in G_K$, let $a=g|_F$ and $b=\phi_m(g)|_E$. The construction of $\theta_E$ is compatible with the corresponding conjugation actions. In particular, by contravariance, $a\tau_E b^{-1}$ induces the same map as $\tau_E$. The faithfulness of $A_2(-)$ (see \cite[Theorem 2.6]{Ho3}) gives $a\tau_E=\tau_E b$.

For $E\subset E'$ in the cofinal family, write $F\subset F'$ for the associated source fields. Naturality of restriction shows that $E\hookrightarrow E'\xrightarrow{\tau_{E'}}F'$ and $E\xrightarrow{\tau_E}F\hookrightarrow F'$ induce the same map $A_2(F')\to A_2(E)$. By faithfulness the two embeddings coincide. Taking the union gives a $\phi_m$-equivariant embedding $\tau_m:L_m\hookrightarrow K_m$, which sends $L$ into $K$ by taking invariants. Conversely, equivariance on the roots of unity gives $\chi_L\circ\phi_{m+2}=\chi_K$.

Put $T=\tau_m(L_m)$ and $k_0=\tau_m(L)$. The compact restriction image of $G_K^m$ in $\operatorname{Gal}(T/k_0)$ fixes precisely $T\cap K$. Thus
\[
 [G_L^m:\phi_m(G_K^m)]=[T\cap K:k_0]\leq[K:k_0].
\]
Any two equivariant embeddings induce the same map $A_2(K)\to A_2(L)$ and hence agree on $L$. Their images of $L_m$ coincide by normality over the common base image. Accordingly they differ by precomposition with an element $z\in G_L^m$ centralising the open subgroup $H=\phi_m(G_K^m)$. If $m\geq2$, for every finite prime $Q$ of $L_m$ one has $H\cap D_Q=H\cap D_{zQ}$. This subgroup is open in both decomposition groups, so $Q=zQ$ by \cite[Proposition 1.9, (iii)$\Rightarrow$(i)]{ST1}. It follows that $z=1$ by \cite[Lemma 1.8]{ST1}. For $m=1$, every $z\in G_L^{\rm ab}$ centralises $H$, giving precisely the stated torsor.
\end{proof}

\begin{cor}\label{cor:isom}
Let $m\geq1$ and let $\phi_{m+2}:G_K^{m+2}\xrightarrow{\sim}G_L^{m+2}$ be an isomorphism. There exists a $\phi_m$-equivariant isomorphism $\sigma_m:L_m\xrightarrow{\sim}K_m$ with $\sigma_m(L)=K$. It is unique for $m\geq2$, and its restriction to $L$ is unique for $m=1$.
\end{cor}
\begin{proof}
By \cite[Theorem 1.26]{ST1}, $\phi_{m+2}$ and its inverse preserve the cyclotomic characters. Apply Theorem~\ref{thm:criterion} in both directions. The composites of the resulting embeddings induce the identity maps on the corresponding $A_2$ groups, so faithfulness gives inverse maps on the base fields. An embedding extending an isomorphism of the base fields identifies their maximal $m$-step solvable extensions, by extending it to an isomorphism of algebraic closures. The uniqueness assertions follow from Theorem~\ref{thm:criterion}.
\end{proof}
\begin{rem}\label{rem:ST-isom}
For every $m\geq1$, Corollary~\ref{cor:isom} lowers the input depth in \cite[Theorem 2]{ST1} from $m+3$ to $m+2$, with the same equivariance and uniqueness assertions. At $m=1$, the corollary also makes \cite[Theorem 1]{ST1} functorial without increasing its input depth. Namely, a given isomorphism $\alpha:G_K^3\xrightarrow{\sim}G_L^3$ determines a unique field isomorphism $\tau_\alpha:L\xrightarrow{\sim}K$ which admits an $\alpha_1$-equivariant extension to the maximal abelian extensions. Although such an extension need not be unique, the composite of two successive equivariant extensions is again equivariant. Thus, for $\beta:G_L^3\xrightarrow{\sim}G_M^3$, uniqueness on the base fields gives $\tau_{\beta\circ\alpha}=\tau_\alpha\circ\tau_\beta$, and the identity induces the identity. The compatibility here is with the induced isomorphism at the abelian layer.
\end{rem}
\begin{rem}
For $m=1$, the construction in the proof of Theorem~\ref{thm:criterion} singles out one embedding from the $G_L^{\rm ab}$-torsor: its restriction to every $E$ in the cofinal family induces the prescribed map $\theta_E$. This condition uses the higher map and is stronger than $\phi_1$-equivariance alone.
\end{rem}
\section{The Logarithmic Cyclotomic Compatibility}
In this section, we prove the logarithmic cyclotomic compatibility in the $m$-step solvable case, we shall follow the strategy in \cite{OAI} and \cite{Tsuji}, with a suitable adaptation in the finite-step solvable situation. \par 
Fix $\ell \in \mathfrak{Primes}$ and a number field $A$, we shall write
$$
X_\ell(A) := \text{Hom}(G_A,\mathbb{Q}_\ell)
$$
for the vector space of all continuous characters $\rho: G_A \to \mathbb{Q}_\ell$. Those characters are unramified at every $\mathfrak{p} \nmid \ell$. If $A/\mathbb{Q}$ is Galois, then the action of $\text{Gal}(A/\mathbb{Q})$ on $X_\ell(A)$ is given by $(gx)(h)=x(\tilde g^{-1}h\tilde g)$ where $\tilde{g} \in G_\mathbb{Q}$ is a lifting of $g$. We shall denote the averaging operator by 
$$
\mathcal{A}_J := |J|^{-1}\sum_{g \in J}g
$$
for any finite group $J$. For a prime number $\ell$, write $\lambda_{A,\ell}$ for the $\ell$-adic logarithmic cyclotomic character, i.e. $\lambda_{A,\ell} = \log_\ell\chi_{A,\ell}$. In particular, we have $X_\ell(\mathbb{Q}) = \mathbb{Q}_\ell\lambda_{\mathbb{Q},\ell}$.
\begin{lem}
Let $K$ be a number field, it holds that
$$
r_2(K)+1\leq\dim X_\ell(K) \leq [K:\mathbb{Q}]
$$
where $r_2(K)$ is the number of conjugate pairs of non-real embeddings of $K$ and the dimension is the $\mathbb{Q}_\ell$-dimension.
\end{lem}
\begin{proof}
See Section 3.1 in \cite{OAI}.
\end{proof}
\begin{lem}
Let $K$ be a number field and let $\beta:G_K\to P^4$ have open image in $P^2$. Fix an odd prime $\ell$. Then all but finitely many rational primes $q\neq\ell$ are targets of primes active at $\ell$.
\end{lem}
\begin{proof}
Suppose that infinitely many rational primes $q_i\neq\ell$ are missed. Put $C_0=\mathbb{Q}(\zeta_\ell)$ and $C_\infty=C_0(q_i^{1/\ell}:i\geq1)\subset\mathbb{Q}^{(2)}$. The valuations at primes over the $q_i$ give independent Kummer classes, so $C_\infty/C_0$ has infinite elementary abelian rank. Let $K'/K$ correspond to the inverse image of the subgroup fixing $C_0$. By openness at depth two, the image of $G_{K'}$ in $\operatorname{Gal}(C_\infty/C_0)$ is open, so the extension of $K'$ defined by this map has infinite elementary abelian rank.

If this extension ramifies at a prime $v'\nmid\ell$ of $K'$, a pro-$\ell$ Sylow subgroup of its local inertia has nontrivial image. Thus $v'$ is active for the restriction of $\beta$. Its inertia maps at depth two into inertia at its target, which must be one of the $q_i$, since $C_\infty/C_0$ is unramified outside $\ell$ and the $q_i$. The prime $v$ below $v'$ in $K$ is also active: after compatible choices, local inertia at $v'$ is the intersection of local inertia at $v$ with the open source subgroup, and the pro-$\ell$ Sylow at $v'$ is contained in one at $v$. Their lower target primes agree by the containment of the nontrivial inertia image in both decomposition groups and \cite[Proposition 1.3]{ST1}. This contradicts the choice of $q_i$. The pulled-back extension is therefore unramified outside $\ell$, contradicting the finiteness of the elementary abelian $\ell$-rank for a fixed finite ramification set.
\end{proof}
We now adapt the method of \cite[Lemma 1.7]{Tsuji} to derive the logarithmic cyclotomic compatibility in the finite-step solvable case.
\begin{thm}
Let $\beta:G_K\to P^4$ have open image in $P^2$. For some integer $n_0\geq1$ and every sufficiently large odd prime $\ell\equiv1\pmod{n_0}$, 
$$
\lambda_{\mathbb{Q},\ell} \circ \beta= \lambda_{K,\ell}. 
$$
\end{thm}
\begin{proof}
Let $M/\mathbb{Q}$ be the Galois closure of $K/\mathbb{Q}$ and put $n_0=\exp\operatorname{Gal}(M/\mathbb{Q})$. Let $A=(\mathbb{Q}^{(2)})^{\beta_2(G_M)}$, which is finite by the assumption on the image. Choose an odd prime $\ell\equiv1\pmod{n_0}$ unramified in $A$. For $b\geq1$, put
\[
 R_b=\mathbb{Q}(\ell^{1/\ell^b}),\qquad
 E_b=R_b\mathbb{Q}(\zeta_{\ell^b}),\qquad R_0=\mathbb{Q}.
\]
The element $\ell$ is not an $\ell$-th power in $\mathbb{Q}(\zeta_{\ell^b})$, since otherwise this abelian field would contain the non-Galois field $\mathbb{Q}(\ell^{1/\ell})$. By Kummer theory $[E_b:\mathbb{Q}(\zeta_{\ell^b})]=\ell^b$, while by Eisenstein $[R_b:\mathbb{Q}]=\ell^b$. Thus the stabiliser of $\ell^{1/\ell^b}$ is a complement and
\[
 \operatorname{Gal}(E_b/\mathbb{Q})\simeq V\rtimes C,\qquad
 V\simeq\mathbb Z/\ell^b\mathbb Z,\quad
 C\simeq(\mathbb Z/\ell^b\mathbb Z)^\times,\quad E_b^C=R_b .
\]
Here $C$ acts on $V$ by multiplication and $c:=|C|=(\ell-1)\ell^{b-1}$. The commutator subgroup is $V$ (the element $2\in C$ already gives all of $V$ by commutation). In particular, $E_b\subset\mathbb{Q}^{(2)}$ and $E_b^{\rm ab}\subset\mathbb{Q}^{(3)}$. By the discriminant of the radical polynomial and the ramification of the cyclotomic field, $E_b$ is unramified outside $\ell$. Thus $A\cap E_b$ is unramified at every finite rational prime and is $\mathbb{Q}$ by the discriminant bound. Consequently $\beta_2(G_M)$ maps onto $\operatorname{Gal}(E_b/\mathbb Q)$ for every $b$.
\par
The field $R_b$ has signature $(1,(\ell^b-1)/2)$, then it follows from Lemma 4.1 that
$$
\dim X_\ell(R_b) \geq (\ell^b+1)/2 > \ell^{b-1} \geq \dim X_\ell(R_{b-1}).
$$
Take $0\neq x \in \ker(\operatorname{Cor}_{R_b/R_{b-1}})$ (with corestriction on the character spaces $X_\ell(-)$). Let $y := \text{Res}_{E_b/R_b}(x)$, the vector $y$ is also non-zero and is fixed by $C$. Then
$$
\text{Res}_{E_b/R_{b-1}}\text{Cor}_{R_b/R_{b-1}}x = \sum_{u \in V[\ell]} uy.
$$
Write $W_0=\mathbb Q_\ell[V]y$. Averaging over $V[\ell]$ acts as the identity on the factors of $\mathbb Q_\ell[V]$ belonging to characters of order less than $\ell^b$ and as zero on the primitive factor $\mathbb Q_\ell(\zeta_{\ell^b})$. The preceding identity therefore places $W_0$ in the primitive part. Since $y\neq 0$ and that factor is a field, $W_0$ is isomorphic to the primitive factor as a $V$-module. In the identification sending $y$ to $1$, $V$ acts by multiplication by roots of unity, and $C$ acts by the corresponding Galois automorphisms. The action of $V\rtimes C$ on $W_0$ is faithful: an element acting trivially must fix $1$, so its $V$-component is trivial, and then its $C$-component is trivial. Also $\dim W_0^C=1$.
\par
Now write $\rho$ for the $G_M$-character $\lambda_{\mathbb{Q},\ell} \circ \beta - \lambda_{M,\ell}$. The values of $\rho$ are in $\mathbb{Z}_\ell$. For $\sigma \in \text{Gal}(M/\mathbb{Q})$, let $(\sigma\cdot\rho)(g)=\rho(\tilde\sigma^{-1}g\tilde\sigma)$ for $g \in G_M$ and any lifting $\tilde{\sigma} \in G_\mathbb{Q}$. Notice that this action does not depend on the choice of the liftings $\tilde{\sigma}$ since $\rho$ is an additive character (which factors through abelianisation). List the distinct characters in this finite orbit as $\rho_1,\dots,\rho_s$. For each $1 \leq i<j\leq s$, choose $g_{ij} \in G_M$ with $\rho_i(g_{ij}) \neq \rho_j(g_{ij})$. Define
$$
b_0 := 1+\max{\{\{0\} \cup \{\text{ord}_\ell(\rho_{i}(g_{ij})-\rho_j(g_{ij})):1 \leq i<j\leq s\}\}}.
$$
For every $b \geq b_0$, the reductions of $\rho_i$ and $\rho_j$ modulo $\ell^b$ differ at $g_{ij}$. In particular, for each $\sigma \in \text{Gal}(M/\mathbb{Q})$, 
$$
\sigma\cdot \rho \equiv \rho \pmod{\ell^b}~\text{on $G_M$} \implies \sigma\cdot \rho = \rho.
$$
We shall fix such an integer $b$. \par 
Choose a finite Galois extension $F/\mathbb{Q}$ containing the fixed field of the kernel of $G_M \twoheadrightarrow V \rtimes C$, in particular $M \subset F$. Write $\Gamma := \text{Gal}(F/\mathbb{Q})$ and $H := \text{Gal}(F/M)$, write $\theta:H \twoheadrightarrow V \rtimes C$ for the induced quotient. Since $E_b^{\rm ab}\subset\mathbb{Q}^{(3)}$, restriction of $\beta|_{G_F}$ defines the pullback
$$
j:X_\ell(E_b) \to X_\ell(F) 
$$
which is $H$-equivariant, with action on the source through $\theta$. We shall prove injectivity on $W_0$ below; no openness at depth three is assumed. \par 
Fix a generator $c_*$ of $C$. By Chebotarev density theorem, there are infinitely many rational primes $q \neq \ell$ whose Frobenius conjugacy class in $\text{Gal}(E_b/\mathbb{Q})$ is the class of $c_*$. For each such $q$, choose $w_q^*$ above $q$ in $E_b$ with Frobenius exactly $c_*$. Since $q$ is unramified, $D_{w_q^*}(E_b/\mathbb{Q}) = C$. Apply Lemma 4.2 to $\beta|_{G_M}$, for all but finitely many such $q$, there is a source prime $v_q$ of $M$ with target rational prime $q$. \par 
For each remaining $q$, choose a prolongation $\tilde{v}_q^0$ of $v_q$ to $\bar{\mathbb{Q}}$, a pro-$\ell$ Sylow subgroup $S_q^0 \subset D_{\tilde{v}_q^0} \subset G_M$. Furthermore choose compatible primes $\mathfrak{m}^0_{3,q},\mathfrak{m}^0_{2,q}:=\mathfrak{m}^0_{3,q} |_{\mathbb{Q}^{(2)}}$ as in Lemma 2.1. Set $w_q^0 := \mathfrak{m}_{2,q}^0|_{E_b}$. It follows from the transitivity of the Galois action on primes, we may choose $h_q \in G_M$ with $\beta_2(h_q)|_{E_{b}}w_q^0 = w_q^*$. Now put
$$
\tilde{v}_q = h_q\tilde{v}_q^0~;~S_q=h_qS_q^0h_q^{-1}~;~\mathfrak{m}_{j,q}=\beta_j(h_q)\mathfrak{m}_{j,q}^0~\text{for}~j=2,3.
$$
As $h_q$ fixes $M$, $\tilde{v}_q$ still lies above $v_q$. Conjugating the two containments in Lemma 2.1 gives
$$
\beta_3(S_q) \subset D_{\mathfrak{m}_{3,q}}(\mathbb{Q}^{(3)}/\mathbb{Q})~\text{and}~\beta_2(D_{\tilde{v}_q}) \subset D_{\mathfrak{m}_{2,q}}(\mathbb{Q}^{(2)}/\mathbb{Q}),
$$
and by construction, we have $\mathfrak{m}_{3,q}|_{\mathbb{Q}^{(2)}} = \mathfrak{m}_{2,q}$ and $w'_q := \mathfrak{m}_{2,q}|_{E_b} = w_q^*$. \par 
Set $v_q' := \tilde{v}_q|_F$ and $p_0(q) := v'_q|_{\mathbb{Q}} = v_q|_{\mathbb{Q}}$. The primes $v_q$ are all distinct since each has a unique target rational prime $q$. Thus $q \mapsto p_0(q)$ has finite fibres. After excluding the $q$ for which $p_0(q)$ ramifies in $F$, choose an infinite set $\mathcal{T}$ of remaining $q$ on which $p_0(q)$ is injective. For $q \in \mathcal{T}$, let $\delta_q := \operatorname{Frob}_{v_q'/p_0(q)} \in \Gamma$ be the Frobenius at the specified prime $v_q'$. Since $\Gamma$ is finite, there exists $\delta \in \Gamma$ such that $\{q \in \mathcal{T}:\delta_q=\delta\}$ is infinite. In this case, we set $\Delta := \langle\delta \rangle$. Also, we write $v,v',w',p_0$ for $v_q,v'_q,w'_q,p_0(q)$ respectively. \par 
Write $f$ for the residue degree of $v$ over $p_0$, equivalently the order of $\delta|_M$. The Frobenius lift at $v$ maps to $\delta^f\in H$. The lower containment gives $\theta(\delta^f)\in C$, and Lemma 2.3 identifies its image under $C\xrightarrow{\sim}(\mathbb Z/\ell^b\mathbb Z)^\times$ with $p_0^f\bmod\ell^b$. Since $f\mid n_0$ and $n_0\mid\ell-1$, we have $f\mid c$. Independently, activity gives $p_0\neq\ell$, so $p_0^c\equiv1\bmod\ell^b$ because $c=|(\mathbb Z/\ell^b\mathbb Z)^\times|$. Therefore $\theta(\delta^c)=\theta(\delta^f)^{c/f}=1$, and $\delta^c$ acts trivially on $j(W_0)$.
\par
For a selected $v'$, let $D_{v'} \subset G_F$ be its decomposition group and $I_{v'}$ for the inertia subgroup for the chosen prolongation $\tilde{v}_q$. And choose a Frobenius lifting $g \in D_{v'}$ of $v'$. In the procyclic group $\overline{\langle g \rangle}$, take its pro-$\ell$ component $g_\ell$. Any additive (logarithmic) $\ell$-adic character takes the same value of $g_\ell$ as on $g$ since the prime-to-$\ell$ part has zero image. We have $g_\ell \in D_{v'} \subset D_v$. Take a pro-$\ell$ Sylow of $D_v$ containing $g_\ell$. By Lemma 2.1, the image of such pro-$\ell$ Sylow is contained a target 3-step solvable decomposition subgroup, whose associated prime restricts to $\mathbb{Q}^{(2)}$ is the one attached to $v$. Thus the image of $g_\ell$ fixes $E_b$, which gives a local element at $w'$ with residue exponent $a \in \widehat{\mathbb{Z}}$. Then Lemma 2.3 gives 
$$
a_\ell \log _\ell Nw' = \log_\ell Nv'
$$
where $a_\ell$ is the $\ell$-adic part of $a$. The logarithms on the both sides are non-zero since $p_0,q\neq \ell$. Then for every $x_0 \in W_0$, we have
$$
j(x_0)(\operatorname{Frob}_{v'}) = \frac{\log_\ell Nv'}{\log_\ell Nw'}x_0(\operatorname{Frob}_{w'}). 
$$
The functional $x_0\mapsto x_0(\operatorname{Frob}_{w'})$ on $W_0$ is $C$-invariant, since $D_{w'}(E_b/\mathbb Q)=C$. Each element of $C$ lifts to the decomposition group over $\mathbb Q$ at the chosen prolongation of $w'$. Conjugation by such a lift is trivial on the unramified quotient, and the characters in $X_\ell(E_b)$ are unramified at $w'$, since $q\neq\ell$. Fix $0\neq y_0\in W_0^C$. It descends to the nonzero character $c^{-1}\operatorname{Cor}_{E_b/E_b^C}y_0$ of $G_{E_b^C}$. The chosen $w'$ restrict to distinct degree-one primes of $E_b^C$, with relative residue degree $c$ in $E_b$. By \cite[Proposition 3.2(i)]{OAI}, $y_0(\operatorname{Frob}_{w'})\neq0$ for all but finitely many of them, since evaluation upstairs multiplies the descended value by $c$. The displayed formula therefore shows $j(y_0)\neq0$. The kernel of $j|_{W_0}$ is $V$-stable, because $H$ maps onto $V\rtimes C$. Irreducibility of $W_0$ under $V$ yields injectivity. Put $W=j(W_0)$.

For $x_0\in W_0$ with $\mathcal A_Cx_0=0$, all the indicated evaluations vanish. The formula and \cite[Proposition 3.2(ii)]{OAI} give $\mathcal A_\Delta j(x_0)=0$, since each $v'$ has decomposition group $\Delta=\langle\delta\rangle$ in $F/\mathbb Q$. Conversely, if $\mathcal A_Cx_0\neq0$, it is a nonzero multiple of $y_0$. Hence $j(x_0)(\operatorname{Frob}_{v'})\neq0$ at all but finitely many of the chosen primes. At each of these primes, $\mathcal A_\Delta j(x_0)$ has the same evaluation as $j(x_0)$, by the action of the decomposition group on the unramified quotient. Thus
\[
 j\bigl(\ker(\mathcal A_C|_{W_0})\bigr)=\ker(\mathcal A_\Delta|_W),
\]
and $\mathcal A_\Delta|_W$ has rank one.
\par
Let $D=[H,H]$. It is normal in $\Gamma$, and $\theta(D)=V$. Put $k_\ell=\mathbb Q_\ell$ and let $\overline{k}_\ell$ be an algebraic closure. Choose a $D$-eigenvector $u$ in $W\otimes_{k_\ell}\overline{k}_\ell$ such that $\mathcal A_\Delta u\neq0$. Since $\mathcal A_\Delta(W)$ is a $k_\ell$-line, rescale $u$ so that $z=\mathcal A_\Delta u$ belongs to $X_\ell(F)$. Put $Z=k_\ell[D]z$. This subspace is stable under $\delta$, since $\delta$ fixes $z$ and normalises $D$.

The span of $\delta^ju$ for $0\leq j<c$ contains $z$ because $\delta^cu=u$. Each $\delta^ju$ is a $D$-eigenvector whose character is obtained from that of $u$ by precomposition with conjugation by $\delta^{-j}$ on $D$. These characters all have order $\ell^b$. Thus the span is $D$-stable, contains $Z\otimes_{k_\ell}\overline{k}_\ell$, and has dimension at most $c$. Since $Z\neq0$, at least one character of order $\ell^b$ occurs in $Z\otimes_{k_\ell}\overline{k}_\ell$ by semisimplicity. Some element of $D$ therefore has a primitive $\ell^b$-th root of unity as an eigenvalue on $Z$. Its minimal polynomial over $k_\ell$ has degree $c$, giving $\dim_{k_\ell}Z\geq c$. Equality forces the two spans to coincide, so $u$ lies in the intersection of the scalar extensions of $Z$ and $W$. Scalar extension commutes with intersections of vector subspaces, so $Z\cap W\neq0$. Irreducibility of $W$ as a $k_\ell[D]$-module and equality of dimensions now give $Z=W$.
\par
The faithfulness of the action of $\text{Gal}(E_b/\mathbb{Q})$ on $W$ shows that $\delta$ normalises $\text{ker}(\theta)$ in $H$ (i.e. $\text{ker}(\theta) = \text{ker}(H \to \text{GL}(W))$). Therefore, it induces an automorphism on $V \rtimes C$. The subgroup $V$ is the derived subgroup and hence is characteristic. The automorphism induced by $\delta$ on the quotient $C$ acts on the image of $C\hookrightarrow\operatorname{Aut}(V)$ by conjugation with the restriction of the automorphism to $V$. As $\operatorname{Aut}(V)$ is abelian and the action of $C$ on $V$ is faithful, the induced automorphism of $C$ is the identity. Thus conjugation by a lifting $\delta$ preserves $\chi_{\mathbb{Q},\ell} \circ \beta \pmod{\ell^b}$ on $G_M$. Thus
$$
(\delta|_M)\cdot \rho \equiv \rho \pmod{\ell^b}
$$
on $G_M$, and by the choice of $b$, we have $(\delta|_M)\cdot \rho = \rho$. \par 
Now $\text{Res}_{F/M}\rho$ is fixed by $\Delta$ and vanishes on the chosen $\operatorname{Frob}_{v'}$ by Lemma 2.3. Then by Proposition 3.2 (ii) in \cite{OAI}, its average is $0$ and since it is invariant, the character itself must be $0$. Since restriction of additive characters to an open subgroup is injective, we obtain the desired equality:
$$
\lambda_{\mathbb{Q},\ell} \circ \beta = \lambda_{K,\ell}.
$$
\end{proof}
\begin{rem}
The proof of Theorem 4.3 uses the auxiliary module of \cite[Lemma 1.7]{Tsuji} in place of the quotient amplification in \cite[Lemma 4.2]{OAI}, without using the extension to the normaliser in \cite[Lemma 4.1]{OAI}. The exceptional coefficient primes are contained in the finite set excluded before the choice of $b$. When the Galois closure $M/\mathbb Q$ is multiquadratic, the congruence condition on $\ell$ is automatic for odd $\ell$.
\end{rem}

\section{From Logarithmic to Full Cyclotomic Compatibility}
Following \cite{OAI}, we deduce full compatibility from Theorem 4.3. This proves Theorem~\ref{intro:automatic}, which, together with Theorem~\ref{intro:criterion}, gives Theorem~\ref{intro:main}. We use the following notation:
\begin{itemize}
    \item Fix a homomorphism $\beta:G_K\to P^4$ with open image in $P^2$. 
    \item Let $\mathcal L$ consist of all sufficiently large odd primes in the progression of Theorem 4.3, with its finite exclusions.
    \item $B := \mathbb{Q}(i)$ and $F_B$ for the finite extension of $K$ obtained by forming the fixed subfield of $\bar{K}$ of the inverse image of $\text{Gal}(\mathbb{Q}^{(4)}/B)$ by $\beta.$ In this case $[F_B:K] \leq 2$. 
    \item $d:= [F_B:\mathbb{Q}]$.
    \item For a given $\ell \in \mathcal{L}$, and an active prime $v$ of $F_B$, let $\mathfrak{q} \in \mathscr{P}_{\mathbb{Q}^{(2)}}^{<\infty}$ be the prime corresponding $v$ and $w := \mathfrak{q}|_B$. In particular, by Lemma 2.1, the image of $D_v$ via $\beta_2$ is contained in $D_\mathfrak{q}(\mathbb{Q}^{(2)}/B)$; and a pro-$\ell$ inertia of $D_v$ is contained in $I_{\mathfrak{q}}(\mathbb{Q}^{(2)}/B)$. In this case, we call $(v,w)$ an active pair at $\ell$, here $w$ is unique by Proposition 1.3 in \cite{ST1}.
    \item For each $w \in \mathscr{P}_B^{<\infty}$, put $V_w:= \{v \in \mathscr{P}_{F_B}^{<\infty}:(v,w)~\text{is active for some }~\ell \in \mathcal{L}\}.$
\end{itemize}
\begin{lem}
For all but finitely many degree $1$ prime $w \in \mathscr{P}_B^{<\infty}$, $V_w$ is finite.    
\end{lem}
\begin{proof}
 We shall adapt the proof of Lemma 5.1 in \cite{OAI} to prove this lemma. Fix $r \in \mathcal{L}$. Since $B^{\text{ab}} \subset \mathbb{Q}^{(2)}$, restricting $\beta$ to $G_{F_B}$ gives the pullback
 $$
 j_B: X_r(B) \to X_r(F_B).
 $$
 The composite
 $$
 c_B := \text{Cor}_{F_B/\mathbb{Q}}\circ j_B:X_r(B) \to X_r(\mathbb{Q})
 $$
 has non-zero kernel since $\dim X_r(B) \geq 2$ and $\dim X_r(\mathbb{Q}) = 1$. Let $0 \neq \psi \in \text{ker}(c_B)$, and set $x_0 := \lambda_{B,r}$. Theorem 4.3 gives 
 $$
 j_Bx_0 = \lambda_{F_B,r}
 $$
 hence $c_B(x_0 ) = d\lambda_{\mathbb{Q},r} \neq 0$. \par 
 It follows from Proposition 3.2 (i) in \cite{OAI} that $\psi(\text{Frob}_w)\neq 0 $ for all but finitely many degree one primes $w \nmid r$. For such $w$ over a rational prime $q$, write $t:= \psi(\text{Frob}_w)/\log_rq \neq 0$. The character $\psi - tx_0$ is unramified at $w$ and vanishes on the decomposition group of $w$. Then it follows from Lemma 2.1 that $j_B(\psi-tx_0)$ vanishes on the decomposition group of each $v\in V_w$ away from $r$. Since
\[
 \operatorname{Cor}_{F_B/\mathbb Q}j_B(\psi-tx_0)=-td\lambda_{\mathbb Q,r}\neq0,
\]
the desired finiteness follows from \cite[Proposition 3.2(iii)]{OAI}.
 \end{proof}
 \begin{lem}
Let $k$ be a number field and let $S \subset \mathscr{P}_k$ be a finite set containing all archimedean primes and primes above $\ell$ for some rational prime $\ell$. Then
$$
H^1(G_k,\mathbb{Q}_\ell(1))_{\text{unramified outside}~S} \xrightarrow{\sim} \mathcal{O}_{k,S}^\times \otimes_{\mathbb{Z}} \mathbb{Q}_\ell.
$$
 \end{lem}
 \begin{proof}
 See Lemma 5.2 in \cite{OAI}.
 \end{proof}
 \begin{prop}
 For every $\ell \in \mathcal{L}$, at most $2d$ primes $w\nmid \ell$ of $B$ fail to appear in an active pair at $\ell$, if including the primes above $\ell$, the number is at most $2d+2$.    
 \end{prop}
 \begin{proof}
Put $F=F_B$ and $H=\operatorname{Gal}(\mathbb{Q}^{(2)}/B)$. Let $J$ be a finite set of missed primes of $B$ away from $\ell$. For each $w\in J$ choose $a_w\in B^\times$ generating a positive power of the ideal $w$. Its Kummer classes at all finite levels are defined on $H$: the extension obtained by adjoining the $\ell^n$-th roots of $a_w$ and its conjugate over $B(\zeta_{\ell^n})$ is Galois over $\mathbb Q$ and abelian over the abelian field $B(\zeta_{\ell^n})$. Their $\mathbb Q_\ell$-classes are linearly independent by the valuations at $J$.

Write $\rho=\chi_{\mathbb Q,\ell}\circ\beta|_{G_F}$. Pullback gives classes in $H^1(G_F,\mathbb Q_\ell(\rho))$ and is injective: restrict first to the open subgroup $U_F=\beta_2(G_F)\subset H$, using restriction--corestriction over $\mathbb Q_\ell$, then inflate along the surjection from $G_F$. By Theorem 4.3, $\epsilon=\rho/\chi_{F,\ell}$ has finite image in $\mu_{\ell-1}$. Let $F'/F$ be the finite cyclic extension cut out by $\epsilon$ and let $\Delta=\operatorname{Gal}(F'/F)$. Restriction to $G_{F'}$ is again injective. On $G_{F'}$ the two coefficient modules agree: $\rho|_{G_{F'}}=\chi_{F',\ell}$ since $G_{F'}=\ker(\epsilon)$. With the ordinary Tate twist on $G_{F'}$, the restricted classes lie in the $\epsilon^{-1}$-eigenspace for $\Delta$, since the twisted action on cohomology is the ordinary action multiplied by $\epsilon$.

Suppose that the class of $a_w$ ramifies at $u\nmid\ell$ in $F'$. At $u$ the coefficient character is unramified, so the class is nonzero on pro-$\ell$ inertia. The prime $v$ below $u$ in $F$ is then active at $\ell$. By Lemma~\ref{lem:local} its pro-$\ell$ inertia maps into inertia at a prime $t\nmid\ell$ of $B$ occurring in an active pair. In particular $t\neq w$. The Kummer class of $a_w$ is unramified at $t$ since $a_w$ is a unit there, a contradiction. Therefore the classes are unramified outside $S'=\{u:u\mid\ell\}\cup\mathscr{P}_{F'}^\infty$.

By Lemma 5.2 these classes lie in $\mathcal O_{F',S'}^\times\otimes_{\mathbb Z}\mathbb Q_\ell$. Equivariant Dirichlet gives a representation with character $\mathbb Q_\ell[S']-\mathbb Q_\ell$; in particular, each $\Delta$-orbit in $S'$ contributes at most one to the multiplicity of $\epsilon^{-1}$ (see \cite[Proposition 5.3]{OAI}). There are at most $d$ archimedean orbits and at most $d$ orbits above $\ell$, since these orbits correspond to places of $F$. Thus $|J|\leq 2d$. As $J$ was arbitrary, at most $2d$ primes away from $\ell$ are missed; at most two more lie over $\ell$.
 \end{proof}
 \begin{lem}
Let $q_1,\dots,q_s$ be distinct rational primes, and for each $j$, let $\mathcal{P}_j$ be a finite set of rational primes not containing $q_j$. For any integers $n_0,d \geq 1$, there are infinitely many odd primes $\ell\equiv 1 \pmod{n_0}$ distinct from every $q_j$ and every prime in $\bigcup_j \mathcal{P}_j$ outside any prescribed finite set of primes such that
$$
p^f \pmod{\ell} \notin \langle q_j \pmod{\ell} \rangle
$$
where $1 \leq j \leq s$, $p \in \mathcal{P}_j$ and $1 \leq f\leq d$.
 \end{lem}
 \begin{proof}
     See Lemma 5.4 in \cite{OAI} for the proof.
 \end{proof}
 \begin{prop}
All but finitely many degree $1$ primes $w$ of $B$ occur in an active pair at some $\ell \in \mathcal{L}$ with a source prime of the same residue characteristic.
 \end{prop}
 \begin{proof}
 We follow \cite[Proposition 5.5]{OAI}. If the assertion fails, discard the exceptional set of Lemma 5.1 and choose $2d+3$ failing primes $w_j$ above distinct rational primes $q_j$. For each $j$, let $\mathcal P_j$ consist of the rational primes below the finitely many possible partners in $V_{w_j}$. Then $q_j\notin\mathcal P_j$. Lemma 5.4 supplies one $\ell\in\mathcal L$ such that every candidate residue norm $p^f$, with $p\in\mathcal P_j$ and $1\leq f\leq d$, lies outside $\langle q_j\bmod\ell\rangle$. For any active pair at $\ell$ with target $w_j$, Lemma 2.3 would put the corresponding $p^f$ in that subgroup after reduction modulo $\ell$. Thus every selected $w_j$ is missed at this $\ell$, contrary to Proposition 5.3.
 \end{proof}
 \begin{prop}
For every $N \geq 1$, put $E_N := B\mathbb{Q}(\zeta_N) \subset \mathbb{Q}^{(1)}$. Let $\pi_N:P^4 \twoheadrightarrow R_N:=\text{Gal}(E_N/\mathbb{Q})$. Then for every $g \in G_K$, we have
$$
\pi_N(\beta(g)) = \pi_N(g|_{\mathbb{Q}^{(4)}}).
$$
 \end{prop}
 \begin{proof}
 We adapt the proof of Proposition 6.1 in \cite{OAI}. Put $\beta_N := \pi_N \circ \beta$. Fix $r \in \mathcal{L}$, choose a finite Galois extension $\tilde{M}/\mathbb{Q}$ containing $F_B,E_N$ and the fixed subfield $\bar{\mathbb{Q}}^{\text{ker}(\beta_N)}$. In particular, $G_{\tilde{M}} \subset G_{F_B}$ and $\beta(G_{\tilde{M}})$ fixes $E_N$. Since $E_N^{\text{ab}} \subset \mathbb{Q}^{(2)}$, there are two injective pullbacks
 $$
 j_\beta,j_{\text{nat}}:X_r(E_N) \hookrightarrow X_r(\tilde{M})
 $$
 where $j_\beta$ is induced by $\beta|_{G_{\tilde M}}$ followed by restriction to $E_N^{\rm ab}$, and $j_{\rm nat}$ by ordinary restriction. The image of $\beta_2(G_{\tilde M})$ is open in $P^2$, since $G_{\tilde M}$ is open in $G_K$ and $\beta_2$ has open image. Its restriction to $\operatorname{Gal}(E_N^{\rm ab}/E_N)$ is therefore open there, and a continuous additive character vanishing on an open subgroup is zero. This proves injectivity of $j_\beta$; for $j_{\rm nat}$ it follows from finiteness of $\tilde M/E_N$. \par 
 Choose infinitely many rational primes $q \neq r$ splitting completely in $\tilde{M}$, omitting the finite exceptional set in Proposition 5.5. For each of such $q$, choose a prime $w$ above $q$ of $B$, and an active pair $(v,w)$ at some prime $\ell \in \mathcal{L}$ with $v$ a prime of $F_B$ above $q$. Let $v'$ be the induced prime of $\tilde{M}$ above $v$ and let $u$ be the induced prime of $E_N$ on the other side. Both $v'$ and $u$ have degree $1$ since $q$ splits completely in $\tilde{M}$. If $g \in D_{v'}$ is a Frobenius lifting, choose a local lift of $\beta_2(g)$ in the corresponding decomposition group. Its residue exponent is initially unknown. We have
 $$
 a_r\log_rq=\lambda_{\mathbb{Q},r}(\beta(g)) = \lambda_{\tilde{M},r}(g) = \log_rq 
 $$
 where $a_r$ is the $r$-adic component of the residue exponent of the chosen lift. Since $\log_rq \neq 0$, $a_r = 1$. \par 
 Since the $\mathbb{Q}_r$-characters are unramified away from $r$, for all $x \in X_r(E_N)$, we have
 $$
 j_{\beta}(x)(\operatorname{Frob}_{v'}) = x(\operatorname{Frob}_u).
 $$
 On the other hand, 
 $$
 j_{\text{nat}}(x)(\operatorname{Frob}_{v'})=x(\text{Frob}_{u_0})
 $$
 where $u_0 = v'|_{E_N}$. Choose $a \in R_N$ with $a(u)=u_0$, the conjugation action (cf. at the beginning of section 4) gives
 $$
 x(\operatorname{Frob}_u) = (ax)(\text{Frob}_{u_0}).
 $$
Since $R_N$ is finite, we may keep infinitely many $q$ with the same $a$. For each $x$, the difference $j_\beta(x) - j_{\text{nat}}(ax)$ vanishes at infinitely many degree $1$ Frobenius elements, so is zero by Proposition 3.2 (i) in \cite{OAI}. Thus,
$$
j_\beta = j_{\text{nat}} \circ a.
$$
Since $\tilde M/\mathbb Q$ is Galois, the $G_K$-action on $X_r(\tilde M)$ is well-defined. Conjugation and pullback give
$$
gj_{\beta}(x ) = j_\beta(\beta_N(g)x),~\text{and}~gj_{\text{nat}}(x) = j_{\text{nat}}(\sigma_N(g)x)
$$
where $\sigma_N $ denotes the composite $G_K \xrightarrow{g \mapsto g|_{\mathbb{Q}^{(4)}}} P^4 \twoheadrightarrow R_N$. Substitute the identification $j_\beta = j_{\text{nat}}\circ a$ gives
$$
a \beta_N(g)x=\sigma_N(g)ax
$$
for all $x$. Since the action of $R_N$ on $X_r(E_N)$ is faithful by Proposition 3.3 in \cite{OAI}, $a\beta_N(g) = \sigma_N(g)a$. On the other hand, $R_N$ is abelian gives $\beta_N(g) = \sigma_N(g)$, as required.
 \end{proof}
We now obtain full cyclotomic compatibility.
\begin{thm}[Theorem~\ref{intro:automatic}]\label{thm:automatic}
Let $\beta:G_K\to P^4$ be a homomorphism whose projection to $P^2$ has open image. Then $\chi_{\mathbb Q}\circ\beta=\chi_K$. In particular, any $\phi_4:G_K^4\to G_L^4$ with open image at depth two preserves the cyclotomic characters.
\end{thm}
\begin{proof}
Applying Proposition 5.6 for all $N\geq1$ proves the first assertion. For the second, compose $G_K\to G_K^4\xrightarrow{\phi_4}G_L^4$ with restriction to $P^4$. For $j=2,4$, the restriction $G_L^j\to P^j$ has image $\operatorname{Gal}(\mathbb Q^{(j)}/(\mathbb Q^{(j)}\cap L))$, an open subgroup, because $L\mathbb Q^{(j)}\subset L_j$. The projection of the composite to $P^2$ therefore has open image. Restriction to $P^4$ preserves the cyclotomic character, so the first assertion applies.
\end{proof}

\begin{thm}[Theorem~\ref{intro:main}]
Let $m\geq2$ and let $\phi_{m+2}:G_K^{m+2}\to G_L^{m+2}$ be a homomorphism. The following are equivalent: $\phi_2$ has open image; $\chi_L\circ\phi_{m+2}=\chi_K$; and $\phi_m$ is induced by a $\phi_m$-equivariant embedding $L_m\hookrightarrow K_m$. This embedding is unique.
\end{thm}
\begin{proof}
Since $m+2\geq4$, openness of $\phi_2$ gives compatibility by applying Theorem~\ref{thm:automatic} to $\phi_4$. Theorem~\ref{thm:criterion} gives the equivalence with the embedding and its uniqueness. Finally, the index formula in that theorem shows that $\phi_m$ is open, so its depth-two projection $\phi_2$ is open.
\end{proof}
\begin{rem}
In the cyclotomic criterion, the lift to depth $m+2$ gives the $(\star_\ell)$-condition at depth $m+1$, and passage to depth $m$ gives the full decomposition-group containment. In the automatic compatibility argument, the auxiliary characters are evaluated at depth three, while the active correspondence and inertia comparison use depth two. This explains the depth required by the present proof; it does not establish optimality of the input depth.
\end{rem}

\begin{eg}[A non-open lift with open image at depth two]\label{eg:nonopen-lift}
For every $m\geq2$, there is an injective endomorphism of $G_{\mathbb Q}^{m+2}$ whose image has infinite index and which induces the identity on $G_{\mathbb Q}^{m+1}$.
\end{eg}
\begin{proof}
Put $n=m+2$, $G=G_{\mathbb Q}^n$, $B=G_{\mathbb Q}^{n-1}$, and $A=\ker(G\to B)=G^{[n-1]}$. Conjugation gives a $B$-action on the abelian group $A$. Let $D$ be the pro-$2$ factor of $B^{[n-2]}$, which is an abelian characteristic subgroup of $B$. Write $A=A_5\times A_{5'}$ for the pro-$5$ and prime-to-$5$ factors. Then
\[
 A_5=A_5^D\oplus M,\qquad M=\overline{[D,A_5]}.
\]
Indeed, on each finite $B$-stable quotient of $A_5$, averaging over the finite image of $D$ gives a projection onto the fixed subgroup with kernel generated by the commutators with $D$. These projections commute with the quotient maps, and are $B$-equivariant since $D$ is normal.

We first show that $M$ has a closed complement in $G$. Let $\pi:E=G/A_{5'}\to B$ and let $T$ be a pro-$2$ Sylow subgroup of $\pi^{-1}(D)$. By the profinite Sylow theorem (see \cite[\S 1.6]{RibesSylow}), $T$ maps onto $D$ and $T\cap A_5=1$, so $\pi^{-1}(D)=A_5T$. Since $\pi^{-1}(D)$ is normal in $E$, every $E$-conjugate of $T$ is conjugate to $T$ inside $A_5T$, hence by an element of $A_5$. It follows that $\pi(N_E(T))=B$. Also $N_E(T)\cap A_5=A_5^D$: a normalising element of $A_5$ has its commutators with $T$ in $A_5\cap T$, so it centralises $T$. If $H$ is the inverse image of $N_E(T)$ in $G$, then $G=AH$ and $A\cap H=A_{5'}\times A_5^D$. Thus $G=M\rtimes H$.

We next show that $[M:5M]=\infty$. Take a finite $2$-group $U$ of derived length $n-1$, and put $C=U^{[n-2]}$. For instance, one may take the upper unitriangular group of degree $2^{n-2}+1$ over $\mathbb F_2$: its $j$-th derived subgroup is generated by the elementary transvections at distance at least $2^j$ above the diagonal, by the elementary commutator formula. Let $V$ be the kernel of the $C$-averaging projection on the regular module $\mathbb F_5[U]$. Normality of $C$ makes $V$ a $U$-module. The action of $C\neq 1$ on the regular basis is nontrivial, so $V\neq0$. For every $v\in V$ one has $v=|C|^{-1}\sum_{t\in C}(v-tv)$, whence $[C,V]=V$.

For every $r\geq 1$, put $W_r=V^{\oplus r}\rtimes U$. Since $[U^{[j]},V]=V$ for $0\leq j\leq n-2$, we have $W_r^{[j]}=V^{\oplus r}\rtimes U^{[j]}$ for $0\leq j\leq n-1$. In particular, $W_r^{[n-1]}=V^{\oplus r}$ and $W_r^{[n]}=1$. By Shafarevich's theorem (see \cite{SchmidtWingberg}), $W_r$ occurs as a Galois group over $\mathbb Q$, so it is a quotient of $G$. The resulting map sends $A_5$ onto $V^{\oplus r}$, and the induced map $B\to U$ sends $D$ onto $C$. Hence $M$ maps onto $[C,V^{\oplus r}]=V^{\oplus r}$. These quotients are killed by $5$ and have unbounded order as $r$ grows, which proves the index assertion.

Now define $\phi_n:G\to G$ by $\phi_n(xh)=(5x)h$ for $x\in M$ and $h\in H$, writing $M$ additively. This is a continuous homomorphism because multiplication by $5$ commutes with the $H$-action. It is injective by \cite[Proposition 1.9]{CS}, which excludes elements of order $5$ in $G$. It induces the identity on $B$ and has image $5M\rtimes H$, of infinite index in $G$, as required.
\end{proof}

\begin{rem}\label{rem:normaliser-obstruction}
The same decomposition $G=M\rtimes H$ gives an automorphism $xh\mapsto(-x)h$ inducing the identity on $B$. If it were conjugation by $g\in G$, the image of $g$ would lie in the centre of $B$, which is trivial by \cite[Corollary 1.7]{ST1}. But then $g\in A$, and $A$ is abelian, so conjugation by $g$ fixes $M$ pointwise. This contradicts the nontrivial action of inversion on the nonzero pro-$5$ group $M$. Thus the assertion that automorphisms of open subgroups are realised by ambient conjugation, used in the proof of \cite[Lemma 4.1]{OAI}, does not carry over to the finite-step quotient even for the whole target. This does not assert that the particular automorphisms arising from a source normaliser fail to be realised.
\end{rem}

\appendix
\section{Logarithmic compatibility at depth three}\label{app:depth-three}
The following argument records what can be deduced without the auxiliary module of Section~4.
\begin{prop}
Let $\beta:G_K\to P^3$ be a homomorphism whose image in $P^2$ is open. For every rational prime $\ell$,
\[
 \operatorname{Cor}_{K/\mathbb Q}(\lambda_{\mathbb Q,\ell}\circ\beta)
   =[K:\mathbb Q]\lambda_{\mathbb Q,\ell}.
\]
In particular, if $\dim_{\mathbb Q_\ell}X_\ell(K)=1$, then $\lambda_{\mathbb Q,\ell}\circ\beta=\lambda_{K,\ell}$.
\end{prop}
\begin{proof}
Put $C_0=\mathbb Q(i,\zeta_\ell)$ and $U=C_0(q^{1/\ell}:q\text{ rational prime},\ q\nmid2\ell)$. The independent valuations at the primes above $q$ show that $U/C_0$ has infinite elementary abelian rank. Moreover $U\subset\mathbb Q^{(2)}$. After restricting to the finite extension $K'/K$ on which $\beta$ fixes $C_0$, its image in $\operatorname{Gal}(U/C_0)$ is open. The pulled-back elementary abelian extension must ramify at infinitely many primes away from $\ell$, since the elementary abelian $\ell$-rank for a fixed finite ramification set is finite.

At each such prime $v'$ of $K'$, the pro-$\ell$ inertia line $J$ has nontrivial image under $\beta_2$. The group $\operatorname{Gal}(\mathbb Q^{(2)}/C_0)$ is torsion-free by \cite[Proposition 1.9]{CS}, since an involution would act nontrivially on $i$. Thus $\beta_2$ is injective on $J$, and on a pro-$\ell$ Sylow subgroup $S$ containing $J$. Lemma~\ref{lem:local}(1) for the map from $G_{K'}$ puts $\beta_2(S)$ as an open subgroup of a pro-$\ell$ Sylow subgroup of a target decomposition group away from $\ell$. The maximal normal abelian subgroup on each side is the corresponding inertia line. Conjugation on this line shows that $\chi_{\mathbb Q,\ell}\circ\beta=\chi_{K',\ell}$ on $S$. For a Frobenius lift at $v'$, its pro-$\ell$ component lies in a choice of such a Sylow subgroup and has the same value as the lift under every additive $\ell$-adic character. It follows that the restriction of $x=\lambda_{\mathbb Q,\ell}\circ\beta-\lambda_{K,\ell}$ vanishes on Frobenius at these $v'$. At the prime $v$ below $v'$, the character $x$ is unramified, and restriction multiplies its Frobenius value by the positive residue degree $f(v'/v)$. Thus $x$ vanishes at infinitely many distinct primes of $K$. By \cite[Proposition 3.2(iii)]{OAI}, $\operatorname{Cor}_{K/\mathbb Q}x=0$. The result follows by restriction--corestriction. When $X_\ell(K)$ is one-dimensional, corestriction is injective because its value on $\lambda_{K,\ell}$ is nonzero.
\end{proof}
\begin{rem}
The conclusion at depth three concerns logarithmic characters only. The primitive summand used in Theorem 4.3 cannot be evaluated on the same metabelian local images: if $N=\operatorname{Gal}(\mathbb Q^{(2)}/E_b)$ and $A_0=(P^2)^{[1]}$, then $V$ acts trivially on $\operatorname{Hom}_{\mathrm{cont}}(N,\mathbb Q_\ell)$. Indeed, $V$ lifts to $A_0$, whose conjugation is trivial on $N\cap A_0$ and on $N/(N\cap A_0)$. Restriction to $N\cap A_0$ gives a two-step filtration with trivial action on both factors, and averaging over $V$ splits it. Thus the nontrivial primitive $V$-module of Section~4 needs characters from the third layer, whereas at input depth three the first containment of Lemma~\ref{lem:local} only controls their second-layer images.
\end{rem}
\bibliography{ref} 
\bibliographystyle{alpha}

\end{document}